\documentclass{article}
\usepackage[english]{babel}
\usepackage[utf8]{inputenc}
\usepackage{amsmath}
\usepackage{amsthm}
\usepackage{graphicx}
\usepackage{amssymb}

\usepackage{color}

\newtheorem{defin}{Definition}
\newtheorem{prop}{Proposition}
\newtheorem{corol}{Corollary}

\newcommand{\E}{\mathbb{E}}
\newcommand{\I}[1]{\mathbb{I}_{\{ #1 \}}}
\newcommand{\Cov}{\mathrm{Cov}}
\newcommand{\bv}{\Big\vert}

\begin{document}

\title
{\bf On 
  Wigner-Ville Spectra \\ and  the Unicity of Time-Varying \\ Quantile-Based Spectral Densities 
}

\author{ Stefan  \textsc{Birr}\thanks{Research supported by  the Sonderforschungsbereich ``Statistical modelling of nonlinear dynamic processes" (SFB~823, Teilprojekt A1, C1) of the Deutsche Forschungsgemeinschaft.} 
 \\
{\small  Ruhr-Universit\" at Bochum}\\
 Holger   \textsc{Dette}\thanks{Research   supported  by the Sonderforschungsbereich ``Statistical modelling of nonlinear dynamic processes" (SFB~823, Teilprojekt A1, C1) of the Deutsche Forschungsgemeinschaft.}
\\
{\small  Ruhr-Universit\" at Bochum}\\
 Marc  \textsc{Hallin}\thanks{Research supported by an Interuniversity Attraction Pole (2012-2017) of the Belgian Science Policy Office.} 
 \\{\small
ECARES,  Universit\' e Libre de Bruxelles } \\
 Tobias {\sc Kley} \thanks{Supported by the Engineering and Physical Sciences Research Council grant no. EP/L014246/1.}  
\\{\small
London School of Economics}\\
and \\
Stanislav {\sc Volgushev}\\
{\small
University of Toronto}}
\date{}

\maketitle

\begin{abstract}
\small{{\bf Abstract:} The unicity of the time-varying quantile-based spectrum proposed in Birr et al.~(2016) is established via an asymptotic representation result involving Wigner-Ville spectra.} \\  \medskip

{\bf Keywords}: \small{Copula-based spectrum, Laplace spectrum, Quantile-based spectrum, Time-varying spectrum, Wigner-Ville spectrum.\\
}
\end{abstract}

\section{Introduction}

Emmanuel Parzen is perhaps best known for  his pioneering contributions to spectral methods in time series analysis and  kernel density estimation---a method which in the signal processing community is often referred to as the {\it Parzen-Rosenblatt window} method. Another lifelong interest of Emmanuel Parzen was quantile-based inference (see, for instance, Parzen~(2004)). In 2009, when  officially retiring from Texas A\&M, he delivered a ``Last Lecture". A Last Lecture is supposed to convey the essential message of a scholarly career, and, quite significantly, Emmanuel Parzen chose to epitomize his sixty-year long activity as a researcher in statistics with the eloquent title  {\it ``Quantiles are Optimal"}. 

When asked to participate in this memorial volume, we thought that there was no better way to pay tribute to Emmanuel Parzen's memory and outstanding achievement than contributing something on quantile-based spectral analysis, at the intersection of his favorite subjects, quantiles and spectral analysis. 

Quantile-based (equivalently,  copula-based) spectral methods recently have attracted renewed attention in the time-series community.   Pioneering work in this direction has been done by Hong~(1999)  and Li~(2008), who coined the name {\it Laplace spectrum}. Similar  ideas were further developed by  Hagemann~(2013), and extended into {\it cross-spectrum} and spectral {\it kernel} concepts by  Dette et al.~(2015) and Kley et al.~(2016), where the asymptotic normality of  smoothed periodogram-based estimators is established.  Along with the stationary marginal distribution, those quantile-based spectra entirely characterize the bivariate distributions of all couples $(X_t, X_{t- h})$, $h=1,2,\ldots$ in a stationary process $\{X_t,\ t\in\mathbb{Z}\}$, hence provide  much more information than classical second-order spectra;  computing them has been made possible  by Kley~(2016) and his R package `quantspec'.  Additional contributions can be found in Li~(2012, 2014), Lee and Rao~(2012) and  Davis et al.~(2013). 

Whether quantile-based or classical, spectral methods typically require long observation periods, and long observation periods,  in general, imply that  stationarity assumptions are violated. 
 In 
Birr et al.~(2016), we therefore introduce a locally stationary version of the concepts and methods initiated in Dette et al.~(2015) and Kley et al.~(2016). The type of local stationarity required in this context  differs substantially from the usual notions proposed, for instance, by   Dahlhaus~(1997), Zhou and Wu~(2009a and b), Roueff and von Sachs~(2011) or Vogt~(2012). Just as the latter,  however, it involves  the choice of families of approximating stationary processes---a choice which, at first sight, could have an impact on the resulting time-varying spectrum. In this paper, we show that this is not the case, and establish the unicity of  the time-varying quantile-based spectra defined there, irrespective of the chosen family of approximating processes.  
 
 The proof of that unicity property relies on an asymptotic representation (Proposition~\ref{Prop}) of the  time-varying quantile-based spectra in terms of appropriate {\it Wigner-Ville spectra} that no longer involve any specific approximating stationary processes---an idea that, mutatis mutandis,  goes back to Theorem~2.2 in Dahlhaus~(1996).
 
 The outline of the paper is as follows. In Section~2, we  briefly introduce  the concept of the Wigner-Ville spectrum. Section~3 similarly presents the locally stationary quantile-based spectra proposed in Birr et al.~(2016). In Section~4, we establish the asymptotic representation result connecting those spectra with the Wigner-Ville ones. The desired unicity follows as a corollary. 
 
 \section{Wigner-Ville spectra.}

The so-called {\it Wigner-Ville spectrum} for harmonizable processes (Martin and Flandrin 1985)  originates in the signal processing literature, and constitutes one of the few links between the engineering  and the statistical  approaches to spectral analysis.   

 Spectral analysis in engineering and physics   usually deals  with (deterministic) signals $s(t)$. The spectrum of a signal, in that context, is  seldom constant over time (think of music or speech),  so that the concept of a time-varying spectrum or a \textit{time-frequency representation} is quite essential. The problem is that no obvious definition for such a concept exists---the situation is somewhat similar to that in probability and mathematical statistics, where   no obvious concept  exists for the spectrum of a nonstationary process. As a result, different approaches have been used. The Wigner-Ville spectrum is based on one of the classical time-frequency representations of a  complex continuous-time   signal~$s(t)$,  called the {\it Wigner distribution},  of the form 
\begin{equation}\label{WVdistrib}
\mathcal{W}(t,\omega) := \frac{1}{2\pi} \int s^\star(t-\frac{1}{2}\tau) s(t+\frac{1}{2}\tau) e^{-i\tau \omega} d\tau,
\end{equation}
where the star indicates complex conjugation. Eugene Wigner's  original motivation (Wigner 1932) was to be able to calculate the quantum correction to the second virial coefficient of a gas, which indicates how it deviates from the {\it ideal gas} law. The Wigner distribution was later used, in connection with characteristic function methods, in the context of signal analysis by Ville~(1948).  One unpleasant feature of $\mathcal{W}(t,\omega) $, however,  is that it can take negative values.

Now, given a nonstationary continuous-time processs $\{X(t)\}$, we can define the covariance kernel as
\[R(s,t) := \text{Cov}(X(s),X(t)). \]
For a nonstationary process, this kernel is no longer a function of the lag~$|t~\!-~\!s|$;  in order to define a local autocovariance function, one can use $R(t-\tau/2,t+\tau/2)$, which measures autocovariance at lag $\tau$ between two process values  centered about a timepoint $t$. A time-varying spectrum then can be obtained as the Fourier transform of $R$ with respect to $\tau$, namely, 
\[ \Psi(t,\omega) := \frac{1}{2\pi}\int R(t-\tau/2,t+\tau/2) e^{-i \tau \omega} d\tau.\]
Spectra of this form have been considered, for example, by Mark~(1970). But  they share the same (main) disadvantage as the Wigner-Ville distribution:  they are not necessarily positive. As a consequence, they  could not compete with the theory of {e\it volutionary spectra} proposed by Priestley~(1965) and the approach was not further pursued.

{Results by Claasen and Mecklenbr\" auger~(1980) and  Flandrin and Escudid~(1980) brought  the attention of researchers back to the Wigner-Ville concepts. They independently proved that the Wigner-Ville distribution indeed enjoys almost all properties that are desirable for a time frequency representation. For a quick overview of these properties, see Hlawatsch and Boudreaux-Bartels~(1992) and, for a systematic approach, Loynes~(1968). Moreover, some of the more important properties are not compatible with the  requirement of a non-negative function, turning what was perceived as the main drawback of the Wigner-Ville spectrum into a necessary evil.  Motivated by these results, Martin~(1982) defines a time-varying spectrum, called the Wigner-Ville spectrum, based on the Wigner-Ville distribution~(\ref{WVdistrib}).}

A process $\{X(t)\}$ is called \textit{harmonizable} if  its autocovariance function can be represented as
\[ R(s,t) = \frac{1}{4\pi^2} \int \int f(\lambda,\mu) e^{i(\lambda s - \mu t)} d\lambda d\mu \]  for some function $f$. Harmonizable processes are a generalization of weakly stationary processes, and were used by Priestley as a starting point for his theory of evolutionary spectra. 
 Martin~(1982) defined the \textit{Wigner-Ville spectrum} of   a harmonizable process $\{X(t)\}$ as 
 \[ W(t,\omega) := \frac{1}{2\pi} \int f(\omega - \tau/2,\omega + \tau/2) e^{i \tau t} d\tau , \]
 justifying his notation by showing that, under appropriate assumptions, 
\[ W(t,\omega) = E\Big[ \int e^{-i\omega \tau} X(t + \tau/2) X^\star(t-\tau/2)d\tau \Big]. \]
Martin and Flandrin~(1985) extended this definition to discrete-time processes  and proposed a class of estimators for it, based on a weighted covariance estimator
\[ \hat R(t+k,t-k) := \sum_{s \in \mathbb{Z}} \Phi(m,2k) X(s + m + k)X^\star(s+m-k),\]
where $\Phi$ is some adequate lag-window function. They derived the first and second moments for those estimators under the assumption of Gaussian processes, but, just as with the theory of evolutionary spectra,  no meaningful asymptotic results could be derived. 

The first theory    allowing for asymptotic results in this nonstationary context was initiated some ten years later by   Dahlhaus~(1996)  with the introduction of {\it locally stationary processes}. Dahlhaus considers  triangular arrays of stationary processes~$\{X_{t,T}\}$ which are ``close" to a family  $\{X_{\theta}(t)\}_{ \theta \in (0,1)}$ whenever $t/T$ is close to $\theta$, and defines the time-varying spectrum $f(u,\omega)$ of   $\{X_{t,T}\}$ at rescaled time point~$u~\!=~\! t/T$ as  the classical spectrum 
\begin{equation} f(u,\omega) := \frac{1}{2\pi} \sum_{k \in \mathbb{Z}} \text{Cov}(X_u(k),X_u(0)) e^{-i \omega k}\label{qspectrum}\end{equation}
of the stationary process $\{X_u(t)\}$.  
Now, as the family of approximating processes $X_u(t)$ is not unique, it is not obvious from the definition that the resulting spectrum is unique. But, considering  the Wigner-Ville spectrum
\[ f_T(u,\omega) := \sum_{s\in\mathbb{Z}} \Cov(X_{[uT-s/2],T},X_{[uT+s/2],T}) e^{-i\omega s} , \]
 Dahlhaus~(1996) proves that,  under appropriate assumptions, 
\[\int_{-\pi}^\pi |f_T(u,\omega) - f(u,\omega)|^2 d\omega = o(1),\]
 as $T\to\infty$, which, provided that $\omega\mapsto f_T(u,\omega)$ and $\omega\mapsto f(u,\omega)$ are continuous,  implies  unicity of $\omega\mapsto f(u,\omega)$.

%

\section{Copula Spectral Densities for locally stationary processes}
 Consider a  triangular array $(X_{t,T},\ 1\leq t\leq T)$, $T\in\mathbb{N}$, of   finite-length realizations   of nonstationary  processes \mbox{$\{X_{t,T},\ t\in\mathbb{Z}\}$}, $T\in\mathbb{N}$. The quantile-based or copula spectral density kernels of a stationary process  are defined  (Dette et al.~(2015);  Kley et al.~(2016)) in terms of its bivariate marginal distribution functions. Therefore, it seems  natural to use  bivariate marginal distribution functions when evaluating, in the definition of local stationarity,  the distance between the nonstationary process $\{X_{t,T}\}$ 
and its stationary approximations~$\{X_{u}(t) \}$. This leads to the following definition.
\begin{defin} \textrm{\em (Birr, Volgushev, Kley, Dette, and Hallin~2016).} 
A triangular array $\{ (X_{t,T})_{t \in \mathbb{Z}} \}_{T \in \mathbb{N}}$ of processes is called {\em locally strictly stationary (of order two)} if there exists a constant $L>0$ and, for every~$u \in (0,1)$,  a strictly stationary process $\{X_{u}(t) , t \in \mathbb{Z}\}$ such that, for every $1\leq r,s \leq T,$
 \begin{equation}\label{eq1}
    \big\Vert F_{r,s;T}(\cdot,\cdot) - G^{u}_{r-s}(\cdot,\cdot)\big\Vert_{\infty} \leq L \Big(
    \max(|r/T-u|,|s/T-u|) + {1}/{T}\Big),
  \end{equation}
where $\Vert  \cdot  \Vert_\infty$ stands for the supremum norm, while $F_{r,s;T}(\cdot,\cdot)$ and $G^{u}_{k}(\cdot,\cdot)$ denote the joint distribution functions of $(X_{r,T}, X_{s,T})$ and $(X_{u}(0),X_{u}(-k))$, respectively.
\end{defin}
With this concept of nonstationarity, we can transfer the stationary quantile-based concepts of Dette et al.~(2015) from $\{X_{u}(t)\}$ to the  nonstationary $\{X_{t,T}\}$,    defining the {\it time-varying copula spectral density}. Define the {\it lag-$h$-copula cross-covariance kernel}  of  $\{X_{u}(t)\}$ as 
\begin{equation}
\gamma^u_h(\tau_1,\tau_2) := \Cov\Big(\I{X_{u}(t) \leq q^u(\tau_1)},\I{X_{u}(t-h)\leq q^u(\tau_2)}\Big) ,\quad \tau_1,\tau_2 \in (0,1),
\end{equation}\label{indic}
where $q^u(\tau)$ denotes $X_{u}(t)$'s marginal quantile of order $\tau$. If we assume that for all $\tau_1,\tau_2,u$ the lag-$h$-covariance kernels $\gamma^u_h(\tau_1,\tau_2)$ are summable, we obtain the {\it time-varying copula spectral density}
\begin{equation} \label{def:tvf}
 \mathfrak{f}^u(\omega,\tau_1,\tau_2) := \frac{1}{2 \pi} \sum_{h = - \infty}^{\infty} \gamma^u_h(\tau_1,\tau_2) e^{-ih\omega},\ \ \tau_1,\tau_2 \in (0,1), \ \ \omega\in(-\pi,\pi ].
\end{equation}
Comparing those definitions with those of the local spectral densities in Dahl\-haus~(1996), we see that the approximating autocovariances appearing there are replaced with  autocovariances involving the copula transform (equivalently, the quantiles) of the approximating processes. This indicates that the local spectral density kernels~(\ref{def:tvf}) can be viewed as  fully non-parametric generalizations of their classical $L^2$-based counterparts,  capturing  pairwise serial dependencies of arbitrary forms. For  detailed comparisons, we refer   to Dette et al.~(2015) and Kley et al.~(2016).

\section{Unicity of the time-varying quantile-based spectral density}
The Wigner-Ville spectrum  associated with   the  series of  indicators involved in the definition (\ref{indic}) of  the copula cross-covariance kernel of ${X_u(t)}$ takes the form 
\begin{equation}\label{WV}\vspace{-10pt} 
\mathfrak{W}_{t_0,T}(\omega,\tau_1,\tau_2 ) = \frac{1 }{2\pi} \sum_{s = - \infty}^{\infty} \gamma_{t_0;T}(s,\tau_1,\tau_2)  e^{-i \omega s},
\end{equation}
where
\begin{equation*}
\gamma_{t_0;T}(s,\tau_1,\tau_2) := \Cov\Big(\I{X_{\lfloor t_0 + s/2\rfloor , T} \leq F^{-1}_{\lfloor t_0 + s/2\rfloor ; T}(\tau_1)}, \I{X_{\lfloor t_0 - s/2\rfloor , T}\leq F^{-1}_{\lfloor t_0 - s/2\rfloor ; T}(\tau_2)}\Big)
\end{equation*}
and $F^{-1}_{t;T}$ denotes the generalized inverse of $F_{t;T}$, the marginal distribution function of $X_{t,T}.$  Indicators being bounded, $\mathfrak{W}_{t_0,T}(\omega,\tau_1,\tau_2 )$  
 exists and is uniquely defined for all $\omega, \tau_1$ and $\tau_2$  as soon as absolute summability holds for the right-hand side of~(\ref{WV}). The following proposition  establishes  a strong relation between the Wigner-Ville spectrum~(\ref{WV}) and the  time-varying copula spectral density kernels~$\mathfrak{f}^u(\omega,\tau_1,\tau_2)
$ defined in~$(\ref{def:tvf})$. To prove that result, we will use the following two regularity assumptions on $X_{t,T}$.

\begin{description}
\item[(A)] There exists $K < \infty$ such that for all $t_0 \in \mathbb{Z}$, $T \in \mathbb{N}$ and $(\tau_1,\tau_2) \in (0,1)^2$, 
\[
\sum_{s\in\mathbb{Z}} |\gamma_{t_0;T}(s,\tau_1,\tau_2)| \leq K.
\]
\item[(C)] For each $t\in \mathbb{Z}$ and $ T\in \mathbb{N}$, the function $x\mapsto F_{t;T}(x)$ is continuous. 
\end{description}
Assumption \textbf{(A)} is a uniform short-range dependence assumption. It holds under typical assumptions on short-range dependency, e.g.\  if the triangular array~$\{X_{t,T}\}$ is $\alpha$-mixing with mixing coefficients $\alpha(k) \leq C k^{-\delta}$ for some $\delta>1$ and~$C$ independent of $T$. Assumption \textbf{(C)} ensures that the marginal distributions $F_{t;T}$ are continuous in a suitable uniform sense.

\begin{prop}\label{Prop} 
Let $\{X_{t,T}\}$ be locally strictly stationary, with approximating processes $\{X_{t}(u)\}$. If moreover $\{X_{t,T}\}$ satisfies \textbf{(A)} and \textbf{(C)}, then
\begin{enumerate}
\item[(a)] the $\gamma_h^u(\tau_1,\tau_2)$'s are absolutely summable for any~$u$ and~$( \tau_1,\tau_2 )\in (0,1)^2$, and 
\item[(b)] for $t_0 = \lfloor uT \rfloor$ and any $(\tau_1,\tau_2) \in (0,1)^2$,
	\[
	\sup_{\omega\in(-\pi,\pi ]} \bv \mathfrak{f}^{u}(\omega,\tau_1,\tau_2 ) -  \mathfrak{W}_{t_0,T}(\omega,\tau_1,\tau_2 ) \bv = o(1)
	\]
as $T\to\infty$.
\end{enumerate}
\end{prop}
\noindent\textbf{Proof.} Throughout this proof, write $t_0 = \lfloor uT \rfloor$. Let $G^u$ and $F_{t,T}$ denote the marginal distribution functions of $X_u(t)$ and $X_{t,T}$ respectively. We begin by proving that under condition \textbf{(C)} the function $x\mapsto G^u(x)$ is continuous for every $u \in (0,1)$. This follows since, for any $x \in \mathbb{R}$ and $ T \in \mathbb{N}$,
\[
\lim_{y \to x} |G^u(x)-G^u(y)| \leq 4L/T  + \lim_{y \to x} |F_{t_0;T}(x)-F_{t_0;T}(y)| = 4L/T.
\] 
This can be made arbitrarily small by choosing $T$ sufficiently large, and continuity of $G^u$ follows. 

Next, observe that expectation of indicators can be written in terms of distribution functions, so that 
\[ \E(\I{X_u(t)\leq x}) = G^u(x) \quad \text{and} \quad \E(\I{X_{u}(t) \leq x}\I{X_{u}(t-h)\leq y}) = G_h^u(x,y), \]
and therefore
\[ \gamma_h^u(\tau_1,\tau_2) = G^{u}_{h}(q^u(\tau_1),q^u(\tau_2)) - \tau_1 \tau_2. \]
Using the same representation for $\gamma_{t_0;T}(h,\tau_1,\tau_2)$, we obtain
\begin{align*}
&|\gamma_h^u(\tau_1,\tau_2)-\gamma_{t_0;T}(h,\tau_1,\tau_2)|\\
& \quad = \Big\vert F_{\lfloor t_0 - h/2\rfloor,\lfloor t_0 + h/2\rfloor;T}(F^{-1}_{\lfloor t_0 - h/2\rfloor;T}(\tau_1),F^{-1}_{\lfloor t_0 + h/2\rfloor;T}(\tau_2)) - G^{u}_{h}(q^u(\tau_1),q^u(\tau_2))\Big\vert .
\end{align*}
Adding and subtracting $ G^{u}_{h}(F^{-1}_{\lfloor t_0 - h/2\rfloor;T}(\tau_1),F^{-1}_{\lfloor t_0 + h/2\rfloor;T}(\tau_2))$,  
and utilizing  the triangle inequality yields 
\begin{flalign*}
&|| F_{\lfloor t_0 - h/2\rfloor,\lfloor t_0 + h/2\rfloor;T}(\cdot,\cdot) - G^{u}_{h}(\cdot,\cdot)||_\infty \\
&\quad + \bv G^{u}_{h}(F^{-1}_{\lfloor t_0 - h/2\rfloor;T}(\tau_1),F^{-1}_{\lfloor t_0 + h/2\rfloor;T}(\tau_2)) - G^{u}_{h}(q^u(\tau_1),q^u(\tau_2))\Big\vert
\end{flalign*}
where  the first term,  in view of local strict stationarity, can be bounded by
\begin{align*}
|| F_{\lfloor t_0 - h/2\rfloor,\lfloor t_0 + h/2\rfloor;T}(\cdot,\cdot) - G^{u}_{h}(\cdot,\cdot)||_\infty & \leq L \Big(\max\big(\big|\frac{t_0+h}{T}-u\big|,\big|\frac{t_0}{T}-u\big|\big)+1/T\Big)\\
& \leq \frac{L(|h| + 2)}{T}.
\end{align*}
For the second term, invoking Sklar's Theorem and the continuity of $G^u$, we can write 
\[G^{u}_{h}(x,y) = \mathcal{C}(G^{u}(x),(G^{u}(y))\]
where $\mathcal{C}$ denotes the copula of $G^{u}_{h}.$ As  copulas are Lipschitz-continuous with constant one, we obtain
\begin{align}\label{RHS}
&\bv G^{u}_{h}(F^{-1}_{\lfloor t_0 - h/2\rfloor;T}(\tau_1),F^{-1}_{\lfloor t_0 + h/2\rfloor;T}(\tau_2)) - G^{u}_{h}(q^u(\tau_1),q^u(\tau_2))\Big\vert\\
&\quad \leq \big|G^u(F^{-1}_{\lfloor t_0 - h/2\rfloor;T}(\tau_1))-G^{u}(q^u(\tau_1))\big|+\big|G^u(F^{-1}_{\lfloor t_0 + h/2\rfloor;T}(\tau_2))-G^{u}(q^u(\tau_2))\big|.
\nonumber\end{align}
  The    continuity of $G^u$ and $F_{\lfloor t_0 - h/2\rfloor;T}$ implies, for the first term in the right-hand side of (\ref{RHS}),   
\begin{align*}
&\big|G^u(F^{-1}_{\lfloor t_0 - h/2\rfloor;T}(\tau_1))-G^{u}(q^u(\tau_1))\big| = \big|G^u(F^{-1}_{\lfloor t_0 - h/2\rfloor;T}(\tau_1))-\tau_1 \big|\\
&\quad = \big|G^u(F^{-1}_{\lfloor t_0 - h/2\rfloor;T}(\tau_1))-F_{\lfloor t_0 - h/2\rfloor;T}(F^{-1}_{\lfloor t_0 - h/2\rfloor;T}(\tau_1))\big|\\
& \quad \leq || F_{\lfloor t_0 - h/2\rfloor;T}(\cdot) - G^{u}(\cdot)||_\infty \leq L(|h|+2)/T,
\end{align*}
where the last inequality holds due to local strict stationarity. Using the same arguments for the second term, we get
\begin{equation}\label{eq:boundgamma}
|\gamma_h^u(\tau_1,\tau_2)-\gamma_{t_0;T}(h,\tau_1,\tau_2)| \leq \frac{3L(|h|+2)}{T}. 
\end{equation}
We first prove part (a) of the proposition. This will be done by contradiction. Assume that, for some $u$, $\tau_1$ and $\tau_2$, the $\gamma_h^u(\tau_1,\tau_2)$'s are not absolutely summable. In this case,   there   exists a $H<\infty$ such that
\begin{equation}\label{eq:low}
\sum_{|h|\leq H} |\gamma_h^u(\tau_1,\tau_2)| > K+1 
\end{equation}
where $K$ is the constant from assumption \textbf{(A)}. On the other hand, given \eqref{eq:boundgamma},
\[
\sum_{|h|\leq H} |\gamma_h^u(\tau_1,\tau_2)| \leq \frac{3L}{T}\sum_{|h|\leq H} (|h|+2) + \sum_{s\in\mathbb{Z}} |\gamma_{t_0;T}(s,\tau_1,\tau_2)| \leq \frac{3L}{T}(2 + 5H + H^2 ) + K.
\] 
By choosing $T$ sufficiently large, this leads to a contradiction with~\eqref{eq:low}. Absolute summability of $\gamma_h^u(\tau_1,\tau_2)$ follows, hence part (a) of Proposition~1. 

We now proceed to prove part~(b). From the absolute summability of the~$\gamma_h^u(\tau_1,\tau_2)$'s, we obtain\vspace{-3mm}
\[
\mathfrak{f}^u(\omega,\tau_1,\tau_2) = \frac{1}{2 \pi} \sum_{h = -T^{1/3}}^{T^{1/3}} \gamma_h^u(\tau_1,\tau_2)e^{-i \omega h} + o(1)
\]
uniformly in $\omega$, while Assumption (A) yields, still uniformly in $\omega$, 
\begin{align*}
\mathfrak{W}_{t_0,T}(\omega,\tau_1,\tau_2 ) &= \frac{1 }{2\pi} \sum_{h = -T^{1/3}}^{T^{1/3}} \gamma_{t_0;T}(h,\tau_1,\tau_2) e^{-i \omega h} + o(1).
\end{align*}
As a consequence of~\eqref{eq:boundgamma}, we have 
\begin{equation*}
\sum_{h = -T^{1/3}}^{T^{1/3}}\left|\gamma_h^u(\tau_1,\tau_2)-\gamma_{t_0;T}(h,\tau_1,\tau_2)\right| = O\left(\frac{T^{2/3}}{T}\right)=o(1), 
\end{equation*}
which establishes the desired result. 
\hfill $\qed$

\medskip

%

The unicity of $
 \mathfrak{f}^u(\omega,\tau_1,\tau_2)$ then follows as an immediate corollary.

\begin{corol}\label{corol} 
For any locally strictly stationary $\{X_{t,T}\}$ fulfilling \textbf{(A)} and \textbf{(C)}, the time-varying copula spectral density $\mathfrak{f}^u(\omega,\tau_1,\tau_2)$ is uniquely defined, i.e.\  does not depend of the choice of the approximating processes $\{X_{t}(u)\}$.
\end{corol}

\end{document}